\documentclass{amsart}
\usepackage{graphicx}
\usepackage{amsfonts}
\usepackage{amscd}
\usepackage{amssymb}
\usepackage{alltt}
\usepackage{multicol}
\usepackage{amsmath}
\usepackage{amsthm}
\usepackage{amscd}

\usepackage[numbers]{natbib}
\usepackage{rotating}
\usepackage{floatpag}
 \rotfloatpagestyle{empty}
\usepackage{amsthm}
\usepackage{graphicx}
\usepackage{multind}
\usepackage{times}

\usepackage{verbatim}
\usepackage{latexsym}
\usepackage{amsfonts}
\usepackage{amsmath}
\usepackage{crop}
\usepackage{fancyhdr}
\usepackage{txfonts}
\usepackage[hyphens]{url}
\usepackage{setspace}
\usepackage{ellipsis} %
\usepackage{wrapfig}
\usepackage{listings}

\usepackage[mathscr]{euscript} 
\usepackage{pifont} 
\usepackage[displaymath]{lineno}

\def\displayallproof{t} 

\def\displayrating{f}

\def\verbose{t}
  {~\par\phantom{!}\endgroup\bigskip}

\def\indy#1#2{\index{index/#1}{#2}\relax}

\def\hide#1{}
\def\swallowed{\relax}
\def\swallow#1\swallowed{}
\newenvironment{iproved}{}{}

\def\hideproof{\renewenvironment{iproved}{%
   \centerline{\it -- Proof Proofed --}
  
  \renewenvironment{enumerate}{}{}
  \def\item{\relax}
  \catcode13=12
  \swallow
}{}}
\def\showproof{\renewenvironment{iproved}{\begin{proof}}{\end{proof}}}
\def\resetproved{\if\displayallproof t\showproof\else\hideproof\fi}

\def\rating#1{\if\displayrating t%
  {{\textsc {[rating={\ensuremath {#1}}].\ }}}\else{}\fi}

\def\oldrating#1{\if\displayrating t%
  {{\textsc {[former rating={\ensuremath {#1}}].\ }}}\else{}\fi}

\def\ifcverbose#1#2{\if\verbose t{{#1}}\else{#2}\fi}

\def\dcg#1#2{{\if\verbose t%
  {{\tt{[DCG-#1]}}\indy{References}{ZC{#2 #1}@{DCG-#1}|page{#2}}}\else{}\fi}}
\def\tlabel#1{\label{#1}\if\verbose t{{\tt [#1].\ }%
   \indy{References}{#1|itt}}\else{}\fi}
\def\formal#1{\relax }

\def\mar#1{}

\def\|{\hbox{\ensuremath{\hspace{0.1em}|\hspace{-0.1em}|\hspace{0.1em}}}}
\def\mid{\ :\ }

\def\normo#1{{\|#1\|}}

\def\=#1{\accent"16 #1}

\newcommand{\ring}[1]{\mathbb{#1}}

\def\v{{\mathbf v}}
\def\u{{\mathbf u}}

\crop

\theoremstyle{plain}

\newtheorem{theorem*}[equation]{Theorem$^*$}

\newtheorem{lemma*}[equation]{Lemma$^*$}

\theoremstyle{definition}

\theoremstyle{remark}

\begin{document}
\title
    {A formal proof of the Kepler conjecture}
\author{Thomas Hales,
  Mark Adams,
  Gertrud Bauer,
  Dang Tat Dat, 
  John Harrison,
  Hoang Le Truong,
  Cezary Kaliszyk,
  Victor Magron, 
  Sean McLaughlin, 
  Nguyen Tat Thang, 
  Nguyen Quang Truong, 
  Tobias Nipkow,
  Steven Obua,
  Joseph Pleso,
  Jason Rute, 
  Alexey Solovyev,
  Ta Thi Hoai An, 
  Tran Nam Trung, 
  Trieu Thi Diep, 
  Josef Urban,
  Vu Khac Ky, 
  Roland Zumkeller,
  }
 
 \lhead{Hales et al.}
\rhead{A formal proof of the Kepler conjecture}



\newcommand{\R}{\mathbb{R}}
\newcommand{\IR}{\mathbb{IR}}
\newcommand{\Arctan}{\mathop{\rm Arctan}}
\newcommand{\abss}[1]{\lvert#1\rvert}
\newcommand{\iabs}{\mathop{\rm iabs}}
\newcommand{\bx}{{\bf x}}
\newcommand{\GT}{G_{\hbox{\tiny Taylor}}}

\begin{abstract}
  This article describes a formal proof of the Kepler conjecture on
  dense sphere packings in a combination of the HOL Light and Isabelle
  proof assistants.  This paper constitutes the official published
  account of the now completed Flyspeck project.
\end{abstract}

\maketitle

\section{Introduction}

The booklet {\it Six-Cornered Snowflake}, which was written by Kepler
in 1611, contains the statement of what is now known as the Kepler
conjecture: no packing of congruent balls in Euclidean three-space has
density greater than that of the face-centered cubic
packing~\cite{Kepler11}.  This conjecture is the oldest problem in
discrete geometry. The Kepler conjecture forms part of Hilbert's 18th
problem, which raises questions about space groups, anisohedral
tilings, and packings in Euclidean space.  Hilbert's questions about
space groups and anisohedral tiles were answered by Bieberbach in 1912
and Reinhardt in 1928. Starting in the 1950s, L.\ Fejes T\'oth gave a
coherent proof strategy for the Kepler conjecture and eventually
suggested that computers might be used to study the
problem~\cite{Fej53}.  The truth of the Kepler conjecture was
established by Ferguson and Hales in 1998, but their proof was not
published in full until 2006~\cite{Hales:2006:DCG}.

The delay in publication was caused by the difficulties that the
referees had in verifying a complex computer proof.  Lagarias has
described the review process~\cite{LagKC}.  He writes, ``The nature of
this proof $\ldots$ makes it hard for humans to check every step
reliably. $\ldots$ [D]etailed checking of many specific assertions
found them to be essentially correct in every case.  The result of the
reviewing process produced in these reviewers a strong degree of
conviction of the essential correctness of this proof approach, and
that the reduction method led to nonlinear programming problems of
tractable size.''  In the end, the proof was published without
complete certification from the referees.

At the Joint Math Meetings in Baltimore in January 2003, Hales announced a
project to give a formal proof of the Kepler conjecture and later
published a project description \cite{hales:DSP:2006:432}.
The project is called {\it Flyspeck}, an expansion of the acronym
FPK, for the Formal Proof of the Kepler conjecture.  The project has
formalized both the traditional text parts of the proof and the
parts of the proof that are implemented in computer code as
calculations.  This paper constitutes the official published account
of the now completed Flyspeck project.  

The first definite contribution to the project was the formal
verification by Bauer and Nipkow of a major piece of computer code
that was used in the proof (see Section~\ref{sec:tc}).  Major work on
the text formalization project started with NSF funding in 2008.  An
international conference on formal proofs and the Flyspeck project
at the Institute of Math (VAST) in Hanoi in 2009 transformed the project
into a large international collaboration.  The PhD theses of Bauer
\cite{Bauer:2006:Thesis}, Obua \cite{Obua:2005:Thesis}, Zumkeller
\cite{roland-thesis}, Magron \cite{Magron:3013:Thesis}, and
Solovyev ~\cite{Solovyev-thesis} have been directly related to this
project.  The book ``Dense Sphere Packings'' gives the mathematical
details of the proof that was formalized~\cite{DSP}.  This article
focuses on the components of the formalization project, rather than
the mathematical details of the proof.

The Flyspeck project has roughly similar size and complexity as other
major formalization projects such as the Feit-Thompson odd order
theorem~\cite{gonthier2013machine}, the CompCert project giving a
verified C compiler \cite{CC}, and the seL4 project that verified an
operating system micro-kernel \cite{Klein-SOSP09}.  The Flyspeck
project might actually set the current record in terms of lines of
code in a verification project.

The code and documentation for the Flyspeck project are available at a
Google code repository devoted to the
project~\cite{website:FlyspeckProject}.  The parts of the project that
have been carried out in Isabelle are available from the Archive of
Formal Proofs~(\url{afp.sf.net}).  Other required software tools are
Subversion (for interactions with the code repository),
Isabelle/HOL~\cite{LNCS2283}, HOL Light~\cite{HOLL}, OCaml (the
implementation language of HOL Light), the CamlP5 preprocessor (for a
syntax extension to OCaml for parsing of mathematical terms), and GLPK
(for linear programming)~\cite{website:GLPK}.


The main statement in the Kepler conjecture can be formally verified
in about five hours on a 2 GHz CPU directly from the proof scripts.
As described in Section~\ref{sec:rec}, the proof scripts of the main
statement have been saved in a recorded proof format.  In this format,
the formal proof of the main statement executes in about forty minutes on a
2 GHz CPU.  We encourage readers of this article to make an
independent verification of this main statement on their computers.


This main statement reduces the Kepler conjecture to three
computationally intensive subclaims that are described in
Section~\ref{sec:statement}.  Two of these subclaims can be checked in
less than one day each.  The third and most difficult of these
subclaims takes about 5000 CPU hours to verify.

\section{The HOL Light proof assistant}\label{sec:hl}

The formal verifications have been carried out in the HOL Light and
Isabelle proof assistants, with a second verification of the main statement
in HOL Zero.

HOL Light is one of a family of proof assistants that implement the
HOL logic \cite{gordon1993introduction}.  It is a classical logic
based on Church's typed lambda calculus and has a simple polymorphic
type system.  For an overview of the deductive system see
\cite{harrison2009hol}.  There are ten primitive inference rules and
three mathematical axioms: the axiom of infinity (positing the
existence of a function $f:X\to X$ that is one-to-one but not onto),
the axiom of extensionality (which states that two functions are equal
if their outputs agree on every input), and the axiom of choice.  HOL
Light has a mechanism that permits the user to extend the system with
new constants and new types.

HOL Light is implemented according to a robust software architecture
known as the ``LCF approach'' \cite{gordon1979edinburgh}, which uses
the type system of its implementation language (OCaml in the case of
HOL Light) to ensure that all deduction must ultimately be
performed by the logic's primitive deductive system implemented in a
kernel.  This design reduces the risk of logical unsoundness to the
risk of an error in the kernel (assuming the correct implementation of
the architecture).

The kernel of HOL Light is remarkably small, amounting to just a few
hundred lines of code.  Furthermore, the kernel has been subject to a
high degree of scrutiny to guarantee that the underlying logic is
consistent and that its implementation in code is free of
bugs~\cite{DBLP:conf/itp/KumarAMO14} \cite{hales-bourbaki2014}.  This
includes a formal verification of its own
correctness~\cite{DBLP:conf/cade/Harrison06}.  It is generally held by
experts that it is extremely unlikely for the HOL Light proof
assistant to create a ``false'' theorem, except in unusual
circumstances such as a user who intentionally hacks the system with
malicious intent.  A formal proof in the HOL Light system is more
reliable by orders of magnitude than proofs published through the
traditional process of peer review.

One feature that makes HOL Light particularly suitable for the
Flyspeck project is its large libraries of formal results for real and
complex analysis.  Libraries include multivariate (gauge) integration,
differential calculus, transcendental functions, and point-set
topology on $\ring{R}^n$.

Proof scripts in the HOL Light are written directly as OCaml code and
are executed in the OCaml read-eval-print loop.  Mathematical
terms are expressed in a special syntax, set apart from ordinary OCaml
code by backquotes.  For example, typed at the OCaml prompt,
\verb!`1`!  equals the successor of zero in the set of natural numbers
(defined from the mathematical axiom of infinity in HOL Light).
It is not to be conflated with the informal $1$ in OCaml, which is
vulnerable to the usual properties of machine arithmetic.

This article displays various terms that are represented in HOL Light
syntax.  For the convenience of the reader, we note some of the
syntactic conventions of HOL Light.  The syntax is expressed entirely
with ASCII characters.  In particular, the universal
quantifier $(\forall)$ is written as (!), the existential quantifier
$(\exists)$ is written (?), the embedding of natural numbers into the
real numbers is denoted (\&), so that for instance \verb!&3! denotes
the real number $3.0$.  The symbol \verb!|-! (the keyboard
approximation for $\vdash$) appears in front of a statement that is a
theorem in the HOL Light system.  Other logical symbols are given by
keyboard approximations: \verb!==>! for $\Longrightarrow$, \verb!/\!
for $\land$, and so forth.

The logic of Isabelle/HOL is similar to that of HOL Light, but it
also includes a number of features not available in HOL Light.  The logic
supports a module system and type classes.  The package includes
extensive front-end support for the user and an intuitive proof
scripting language.  Isabelle/HOL supports a form of computational
reflection (which is used in the Flyspeck project) that allows
executable terms to be exported as ML and executed, with the results
of the computation re-integrated in the proof assistant as theorems.

\section{The statement}\label{sec:statement}

As mentioned in the introduction, the Kepler conjecture asserts that
no packing of congruent balls in Euclidean three-space can have
density exceeding that of the face-centered cubic packing.  That
density is $\pi/\sqrt{18}$, or approximately $0.74$.  The
hexagonal-close packing and various packings combining layers from the
hexagonal-close packing and the face-centered cubic packing also
achieve this bound.  The theorem that has been formalized does not
make any uniqueness claims.

The density of a packing is defined as a limit of the density obtained
within finite containers, as the size of the container tends to
infinity.  To make the statement in the formal proof as simple as
possible, we formalize a statement about the density of a packing
inside a finite spherical container.  This statement contains an error
term.  The ratio of the error term to the volume of the container
tends to zero as the volume of the container tends to infinity.  Thus
in the limit, we obtain the Kepler conjecture in its traditional form.

As a ratio of volumes, the density of a packing is scale invariant.
There is no loss of generality in assuming that the balls in the
packing are normalized to have unit length.  We identify a packing of
balls in $\ring{R}^3$ with the set $V$ of centers of the balls, so
that the distance between distinct elements of $V$ is at least $2$,
the diameter of a ball.  More formally, we have the following
theorem that characterizes a packing.

\begin{obeylines}

\begin{verbatim}

|- packing V <=> 
     (!u v. u IN V /\ v IN V /\ dist(u,v) < &2 ==> u = v)

\end{verbatim}
\end{obeylines}
This states that $V\subset \ring{R}^3$ is a packing if and
only if for every $\u, \v \in V$, if the distance from $\u$ to $\v$ is
less than $2$, then $\u=\v$.  To fix the meaning of what is to be
formalized, we define the constant {\tt the\_kepler\_conjecture} as
follows:
\begin{obeylines}

\begin{verbatim}
|- the_kepler_conjecture <=>
     (!V. packing V
       ==> (?c. !r. &1 <= r
           ==> &(CARD(V INTER ball(vec 0,r))) <=
               pi * r pow 3 / sqrt(&18) + c * r pow 2))
\end{verbatim}
\end{obeylines}
\noindent
In words, we define the Kepler conjecture to be the following claim:
for every packing $V$, there exists a real number $c$ such that for
every real number $r\ge 1$, the number of elements of $V$ contained in an
open spherical container of radius $r$ centered at the origin is at
most
\[
  \frac{\pi\, r^3}{\sqrt{18}} + c\, r^2.
\]
An analysis of the proof shows that there exists a
small computable constant $c$ that works uniformly for all packings
$V$, but we only formalize the weaker statement that allows $c$ to
depend on $V$.  The restriction $r\ge 1$, which bounds $r$ away from
$0$, is needed because there can be arbitrarily small containers whose
intersection with $V$ is nonempty.

The proof of the Kepler conjecture relies on a combination of
traditional mathematical argument and three separate bodies of
computer calculations.  The results of the computer calculations have
been expressed in precise mathematical terms and specified formally in
HOL Light.  The computer calculations are as follows.
\begin{enumerate}
\item The proof of the Kepler conjecture relies on nearly a thousand
  nonlinear inequalities.  
  The term \verb!the_nonlinear_inequalities! in HOL Light is the
  conjunction of these nonlinear inequalities.  See
  Section~\ref{sec:ni}.
\item The combinatorial structure of each possible counterexample to
  the Kepler conjecture is encoded as a plane graph satisfying a
  number of restrictive conditions.  Any graph satisfying these
  conditions is said to be {\it tame}.  A list of all tame plane graphs up
  to isomorphism has been generated by an exhaustive computer search.
  The formal statement that every tame plane graph is isomorphic to one of
  these cases can be expressed in HOL Light as
  \verb!import_tame_classification!.  See Section~\ref{sec:tc}.
\item The final body of computer code is a large collection of linear
  programs.  The results have been formally specified as
  \verb!linear_programming_results! in HOL Light.  See Section~\ref{sec:lp}.
\end{enumerate}

It is then natural to break the formal proof of the Kepler conjecture
into four parts: the formalization of the text part (that is, the
traditional non-computer portions of the proof), and the three
separate bodies of computer calculations.  Because of the size of the
formal proof, the full proof of the Kepler conjecture has not been
obtained in a single session of HOL Light.  What we formalize in a
single session is a theorem

\begin{obeylines}

\begin{verbatim}
|-  the_nonlinear_inequalities /\
    import_tame_classification
    ==> the_kepler_conjecture
\end{verbatim}

\end{obeylines}

This theorem represents the formalization of two of the four parts of
the proof: the text part of the proof and the linear programming.  It
leaves the other two parts (nonlinear inequalities and tame
classification) as assumptions.  The formal proof of the assumption
\verb!the_nonlinear_inequalities! is described in
Section~\ref{sec:ni}.  The formal proof of
\verb!import_tame_classification! in Isabelle is described in
Section~\ref{sec:tc}.  Thus, combining all these results from various
sessions of HOL Light and Isabelle, we have obtained a formalization
of every part of the proof of the Kepler conjecture.

\section{Text formalization}\label{sec:tf}

The next sections of this article turn to each of the four parts of
the proof of the Kepler conjecture, starting with the text part of the
formalization in this section.  In the remainder of this article, we
will call the proof of the Kepler conjecture as it appears in
\cite{Hales:2006:DCG} the {\it original proof}.
The proof that appears in \cite{DSP} will be called the {\it blueprint
  proof}.  The formalization closely follows the blueprint proof.  In
fact, the blueprint and the formalization were developed together, with
numerous revisions of the text in response to issues that arose in the
formalization.

Since the blueprint and formal proof were developed at the same time,
the numbering of lemmas and theorems continued to change as project
took shape.  The blueprint and formal proof are cross-linked by a
system of identifiers that persist through project revisions, as
described in \cite[\textsection 5]{FlyspeckWiki}.

\subsection{standard libraries}

The Flyspeck proof rests on a significant body of mathematical
prerequisites, including basics of measure, geometric notions, and
properties of polyhedra and other affine and convex sets in
$\ring{R}^3$. HOL Light's standard library, largely under the
influence of Flyspeck, has grown to include a fairly extensive theory
of topology, geometry, convexity, measure and integration in
$\ring{R}^n$. Among the pre-proved theorems, for example, are the
Brouwer fixed-point theorem, the Krein-Milman theorem and the
Stone-Weierstrass theorem, as well as a panoply of more forgettable
but practically indispensable lemmas.

Besides the pre-proved theorems, the library includes a number of
automated procedures that can make formal proofs much less painful. In
particular, one tool supports the common style of picking convenient
coordinate axes `without loss of generality' by exploiting
translation, scaling and orthogonal transformation. It works by
automatically using a database of theorems asserting invariance of
various properties under such transformations
\cite{harrison-wlog}. This is invaluable for more intricate results in
geometry, where a convenient choice of frame can make the eventual
algebraic form of a geometrical problem dramatically easier.

The text formalization also includes more specialized results such as
measures of various shapes that are particularly significant in the
partitioning of space in the Flyspeck proof. Among these, for example,
is the usual `angle sum' formula for the area of a spherical triangle,
or more precisely, the volume of the intersection of its conic hull
and a ball. Flyspeck also uses results about special combinatorial and
geometrical objects such as {\em hypermaps} and {\em fans}, and a
substantial number of non-trivial results are proved about these,
including in effect the Euler characteristic formula for planar
graphs.

\subsection{blueprint proof outline}

It is not our purpose to go into the details of the proof of the
Kepler conjecture, for which we refer the reader to \cite{DSP}.  We
simply recall enough of the high level strategy to permit an
understanding of the structure of the formalization.  The proof
considers an arbitrary packing of congruent balls of radius $1$ in
Euclidean space and the properties that it must have to be a
counterexample to the Kepler conjecture.  Recall that we identify a
packing of congruent balls with the discrete set $V$ of centers of the
balls.

The first step of the proof reduces the problem from a packing $V$
involving a countably infinite number of congruent balls to a packing
involving at most finitely many balls in close proximity to one further
fixed central ball.  This reduction is obtained by making a geometric
partition of Euclidean space that is adapted to the given packing.
The pieces of this partition are called Marchal cells (replacing the
decomposition stars in the original proof)~\cite{Marchal11}.  A key
inequality (the Cell Cluster Inequality, which is discussed in
Section~\ref{sec:ab}) provides the link between the Kepler conjecture
and a local optimization problem (the local annulus inequality)
involving at most finitely many balls.  This reduction appears as the
following formal result:


\begin{verbatim}
|- the_nonlinear_inequalities /\ 
   (!V. cell_cluster_inequality V) /\
   (!V. packing V /\ V SUBSET ball_annulus 
        ==> local_annulus_inequality V) 
   ==> the_kepler_conjecture
\end{verbatim}
All three assumptions can be viewed as explicit optimization problems
of continuous functions on compact spaces.
The constant {\it ball annulus} is defined as the set
$A=\{\bx\in\ring{R}^3 \mid 2\le \normo{\bx} \le 2.52\}$.  As this set is
compact and any packing $V$ is discrete, the intersection of a packing
with this set is necessarily finite.   The local annulus inequality for a finite
packing $V\subset A$ can be stated in the following form:
\begin{equation}\label{eqn:la}
\sum_{\v\in V} f(\normo{\v})  \le 12,
\end{equation}
where $f(t) = (2.52-t)/(2.52-2)$ is the linear function that decays
from $1$ to $0$ on the given annulus.  An argument which we do not
repeat here shows that it is enough to consider $V$ with at most $15$
elements~\cite[Lemma~6.110]{DSP}.  We discuss the nonlinear
inequalities and the cell cluster inequality elsewhere in this paper.
The local annulus inequality has been used to solve other open
problems in discrete geometry: Bezdek's strong dodecahedral conjecture
and Fejes-T\'oth's contact conjecture~\cite{hales2011strong}.

The rest of the proof consists in proving the local annulus
inequality.  To carry this out, we assume that we have a
counterexample $V$.  The compactness of the ball annulus allows us to
assume that the counterexample has a special form that we call {\it
  contravening}.  We make a detailed study of the properties of a
contravening $V$.  The most important of these properties is expressed
by what is called the main estimate (see Section~\ref{sec:ab}).  These
properties imply that $V$ gives rise to a combinatorial structure
called a {\it tame planar hypermap}.  (In this brief summary, we
consider tame planar hypermaps to be essentially the same as plane
graphs with certain restrictive properties.  The nodes of the graph
are in bijection with $V$.)  A computer is used to enumerate all of
the finitely many tame plane graphs up to isomorphism.  This reduces
the possible counterexamples $V$ to an explicit finite family of
cases.

For each explicitly given tame planar hypermap $H$, we may consider
all contravening packings $V$ associated with $H$.  By definition,
these are counterexamples to Equation~(\ref{eqn:la}).  We express the
conditions on $V$ as a system of nonlinear inequalities.  We relax the
nonlinear system to a linear system and use linear programming
techniques to show that the linear programs are infeasible and hence
that the nonlinear system is inconsistent, so that the potential
counterexample $V$ cannot exist.  Eliminating all possible
counterexamples in this way, we conclude that the Kepler conjecture
must hold.

\subsection{differences between the original proof and the blueprint
  proof}

The blueprint proof follows the same general outline as the original
proof.  However, many changes have been made to make it more suitable
for formalization.  We list some of the primary differences between
the two proofs

\begin{enumerate}
\item In the blueprint proof, topological results concerning plane
  graphs are replaced with purely combinatorial results about
  hypermaps.  
\item The blueprint proof is based on a different geometric partition
  of space than that originally used.  Marchal introduced this
  partition and first observed its relevance to the Kepler
  conjecture~\cite{Marchal11}.  Marchal's partition is described by
  rules that are better adapted to formalization than the original.
 \item In a formal proof, every new concept comes at a cost: libraries
   of lemmas must be developed to support each concept.  We have
   organized the blueprint proof around a small number of major
   concepts such as spherical trigonometry, volume, hypermap, fan,
   polyhedra, Voronoi partitions, linear programming, and nonlinear
   inequalities of trigonometric functions.
 \item The statements of the blueprint proof are more precise and
   make all hypotheses explicit.
 \item To permit a large collaboration, the chapters of the
   blueprint have been made as independent from one another as
   possible, and long proofs have been broken up into series of
   shorter lemmas.
\item In
  the original, computer calculations were a last resort after
  as much was done by hand as feasible.
  In the
  blueprint, the use of computer has been fully embraced.  As a
  result, many laborious lemmas of the original proof can be
  automated or eliminated altogether.
\end{enumerate}

Because the original proof was not used for the formalization, we
cannot assert that the original proof has been formally verified to be
error free.  Similarly, we cannot assert that the computer code for
the original proof is free of bugs.  The detection and correction of
small errors is a routine part of any formalization project.  Overall,
hundreds of small errors in the proof of the Kepler conjecture were
corrected during formalization.  

\subsection{appendix to the blueprint}\label{sec:ab}

As part of the Flyspeck project, the blueprint proof has been
supplemented with an $84$ page unpublished appendix filled with
additional details about the proof.   We briefly describe the
appendix.

The first part of the appendix gives details about how to formalize
the main estimate~\cite[Sec~7.4]{DSP}.  The main estimate is the most
technically challenging part of the original text, and its
formalization was the most technically challenging part of the
blueprint text.  In fact, the original proof of the main estimate
contained an error that is described and corrected in \cite{HalesHMNOZ10}.
This is the most significant error that was found.

A second major part of the appendix is devoted to the proof of the
Cell Cluster Inequality~\cite[Thm~6.93]{DSP}.   In the blueprint text, the proof is
skipped, with the following comment: ``The proof of this cell cluster
inequality is a computer calculation, which is the most delicate
computer estimate in the book.  It reduces the cell cluster inequality
to hundreds of nonlinear inequalities in at most six variables.''
The appendix gives full details.

A final part of the appendix deals with the final integration of the
various libraries in the project.  As a large collaborative effort
that used two different proof assistants, there were many small
differences between libraries that had to be reconciled to obtain a
clean final theorem.  The most significant difference to reconcile was
the notion of planarity as used in classification of tame plane graphs
and the notion of planarity that appears in the hypermap library.
Until now, in our discussion, we have treated tame plane graphs and
tame planar hypermaps as if there were essentially the same.  However, in
fact, they are based on two very different notions of planarity.  The
recursive procedure defining tame graph planarity is discussed in
Section~\ref{sec:tg}.  By contrast, in the hypermap library, the Euler
characteristic formula is used as the basis of the definition of
hypermap planarity.  The appendix describes how to relate the two
notions.  This part of the appendix can be viewed as an expanded
version of Section~4.7.4 of \cite{DSP}.

\subsection{recording and replaying the proof}\label{sec:rec}

We also have an alternative weaker form of the main statement of the
Kepler conjecture (from Section \ref{sec:statement}) that makes
explicit all three computational portions of the proof.  This form of
the theorem leaves the list of graphs in the archive as an unspecified
bound variable \verb!a!.  This weaker form has the advantage that it
can be verified relatively quickly.

\begin{obeylines}

\begin{verbatim}
|- !a. tame_classification a /\
    good_linear_programming_results a /\ 
    the_nonlinear_inequalities
    ==> the_kepler_conjecture
\end{verbatim}

\end{obeylines}

We have employed an adapted version of HOL Light to record and export
the formal proof steps generated by the proof scripts of the main
statement in this form (\url{www.proof-technologies.com/flyspeck/}).
The exported proof objects have been imported and executed in a
separate HOL Light session involving just the HOL Light core system
and a simple proof importing mechanism.  They have also been imported
and replayed in HOL Zero, another member of the HOL
family~\cite{adams2010introducing}.

This exercise has various benefits.  First, it is generally much
faster to import and replay a recorded formal proof than to execute
the corresponding high-level proof script, whose execution typically
involves a substantial amount of proof search beyond the actual formal
proof steps. Second, it gives a further check of the formal proof.
The successful import into HOL
Zero, a system that pays particular attention to trustworthiness,
effectively eliminates any risk of error in the proof.
Finally, the exported proof objects are available to be imported and
replayed in other HOL systems, opening the Flyspeck libraries of
theorems to the users of other proof assistants.


\section{Nonlinear inequalities}\label{sec:ni}

All but a few nonlinear inequalities in the Flyspeck project have the
following general form
\begin{equation}\label{eqn:D}
\forall \bx,\ \bx \in D \implies f_1(\bx) < 0 \vee \ldots \vee f_k(\bx) < 0.
\end{equation}
where $D = [a_1,b_1] \times \ldots \times [a_n, b_n]$ is a rectangular
domain inside $\R^n$ and $\bx = (x_1,\ldots,x_n)$. In the remaining
few inequalities, the form is similar, except that $k=1$ and the
inequality is not strict.  In every case, the number of variables $n$
is at most six. The following functions and operations appear in
inequalities: basic arithmetic operations, square root, sine, cosine,
arctangent, arcsine, arccosine, and the analytic continuation of
$\arctan(\sqrt{x})/\sqrt{x}$ to the region $x > -1$.  For every point
$\bx\in D$, at least one of the functions $f_i$ is analytic in a
neighborhood of $\bx$ (and takes a negative value), but situations
arise in which not every function is analytic or even defined at every
point in the domain.

The formal verification of inequalities is based on interval
arithmetic~\cite{moore2009introduction}.  For example, $[3.14, 3.15]$
is an interval approximation of $\pi$ since $3.14 \le \pi \le 3.15$.
In order to work with interval approximations, arithmetic operations
are defined over intervals. Denote the set of all intervals over $\R$
as $\IR$. A function (operation) $F:\IR \to \IR$ is called an interval
extension of $f:\R \to \R$ if the following condition holds
\begin{equation*}
\forall I \in \IR,\ \{ f(x)\,:\,x \in I \} \subset F(I).
\end{equation*}
This definition can be easily extended to functions on $\R^k$. It is
easy to construct interval extensions of basic arithmetic operations
and elementary functions. For instance,
\begin{equation*}
[a_1,b_1] \oplus [a_2,b_2]= [a, b] \text { for some $a \le a_1 + a_2$ and $b \ge b_1 + b_2$}.
\end{equation*}
Here, $\oplus$ denotes an interval extension of $+$. We do not define
the result of $\oplus$ as $[a_1 + a_2, b_1 + b_2]$ since we may want
to represent all intervals with limited precision numbers (for
example, decimal numbers with at most 3 significant figures). With
basic interval operations, it is possible to construct an interval
extension of an arbitrary arithmetic expression by replacing all
elementary operations with corresponding interval extensions. Such an
interval extension is called the natural interval extension. Natural
interval extensions can be imprecise, and there are several ways to
improve them. 

One simple way to improve upon the natural interval extension of a
function is to subdivide the original interval into subintervals and
evaluate an interval extension of the function on all subintervals.
Using basic interval arithmetic and subdivisions, it would
theoretically be possible to prove all nonlinear inequalities that
arise in the project. This method does not work well in practice
because the number of subdivisions required to establish some
inequalities would be enormous, especially for multivariate
inequalities.

Both the  C++ informal verification code (from
the original proof of the Kepler conjecture) and our formal
verification procedure implemented in OCaml and HOL Light use improved
interval extensions based on Taylor approximations.

Suppose that a function $g:\ring{R}\to\ring{R}$ is twice differentiable. Fix $y \in
[a,b]$. Then we have the following formula for all $x \in [a,b]$:
\begin{equation*}
g(x) = g(y) + g'(y)(x - y) + \frac{1}{2}g''(\xi)(x - y)^2
\end{equation*}
for some value $\xi = \xi(x) \in [a,b]$.  Choose $w$ such that $w \ge
\max\{y - a, b - y\}$ and define $r(x) = \abss{g'(y)}w +
\frac{1}{2}\abss{g''(\xi(x))}w^2$. We get the following inequalities
\begin{equation*}
\forall x \in [a,b].\quad\ g(y) - r(x) \le g(x) \le g(y) + r(x).
\end{equation*}
Let $G$, $G_1$, and $G_2$ be any interval extensions of $g$, $g'$, and
$g''$. Choose $e$ such that $e \ge \iabs\bigl(G_1([y,y])\bigr)w +
\frac{w^2}{2}\iabs\bigl(G_2([a,b])\bigr)$, where
$\iabs\bigl([c,d]\bigr) = \max\{\abss{c}, \abss{d}\}$.  Assume that
$G([y,y]) = [g_l,g_u]$. Then the following function defines a (second
order) Taylor interval approximation of $g(x)$:
\begin{equation*}
\GT([a,b]) = [l, u] \text{ where $l \le g_l - e$ and $u \ge g_u + e$}.
\end{equation*}
That is,
\begin{equation*}
\forall x \in [a,b].\quad \quad g(x) \in \GT([a,b]).
\end{equation*}
In our verification procedure, we always take $y$ close to the
midpoint $(a+b)/2$ of the interval in order to minimize the value of
$w$. There is an analogous Taylor approximation for multivariate
functions based on the multivariate Taylor theorem.

As a small example, we compute a Taylor interval approximation of
$g(x) = x - \arctan(x)$ on $[1,2]$. We have $g'(x) = 1 - {1}/{(1 + x^2)}$ and $g''(x)
= {-2x}/{(1 + x^2)^2}$. Take $y = 1.5$ and natural interval
extensions $G_1$ and $G_2$ of $g'(x)$ and $g''(x)$. Then $w = 0.5$,
$G([1.5,1.5]) = [0.517,0.518]$, $G_1([1.5,1.5]) = [0.692,0.693]$, and
$G_2([1,2]) = [-0.5, -0.16]$. We get $e = 0.409$ and hence $\GT([1,2])
= [0.108, 0.927]$. This result is much better than the result obtained
with the natural extension $G([1,2]) = [-0.11, 1.22]$. 
We note that in the calculation of $\GT([1,2])$, we evaluate the
expensive interval extension of $\arctan$ only once.  To obtain
similar accuracy with subdivision, more than one evaluation of $\arctan$
is needed.  In general, it is necessary to subdivide the original
interval into subintervals even when Taylor interval approximations are
used. But in most cases, Taylor interval approximations lead to fewer
subdivisions than natural interval extensions.

Taylor interval approximations may also be used to
prove the monotonicity of functions.  By expanding $g'(x)$ in a 
Taylor series, we obtain a Taylor interval approximation in a similar way:
\begin{equation*}
\forall x \in [a,b]. \quad g'(x) \in \GT'([a,b])
\end{equation*}
for some interval $\GT'([a,b])$.  If $0$ is not in this interval, then the
derivative has fixed sign, and
the function $g$ is monotonic, so that the maximum value of $g$ occurs
at the appropriate endpoint.  More generally, in multivariate inequalities,
a partial derivative of fixed sign may be used to reduce the verification on a
rectangle of dimension $k$ to an abutting rectangle of dimension $k-1$.

A few of the inequalities are sharp.  That is, the inequalities to be
proved have the form $f \le 0$, where $f(\bx_0)=0$ at some point $\bx_0$
in the domain $D$ of the inequality.  In each case that arises, $\bx_0$ lies at a
corner of the rectangular domain.  We are able to prove these
inequalities by showing that (1) $f(\bx_0) = 0$ by making a direct
computation using exact arithmetic; (2) $f < 0$ on the complement of
some small neighborhood $U$ of $\bx_0$ in the domain; and (3) every
partial derivative of $f$ on $U$ has the appropriate sign to make the
maximum of $f$ on $U$ to occur at $\bx_0$.  The final two steps are
carried out using our standard tools of Taylor intervals.

All the ideas presented in the discussion above have been
formalized in HOL Light and a special automatic verification procedure
has been written in the combination of OCaml and HOL Light for
verification of general multivariate nonlinear inequalities. This
procedure consists of two parts. The first part is an informal search
procedure which finds an appropriate subdivision of the original
inequality domain and other information which can help in the formal
verification step (such as whether or not to apply the monotonicity
argument, which function from a disjunction should be verified,
etc.). The second part is a formal verification procedure which takes
the result of the informal search procedure as input and produces a
final HOL Light theorem. A detailed description of the search and
verification procedures can be found
in~\cite{Solovyev-thesis,Solovyev:NFM2013}.

All formal numerical computations are done with special finite
precision floating-point numbers formalized in HOL Light. It is
possible to change the precision of all computations dynamically, and
the informal search procedure tries to find minimal precision
necessary for the formal verification. At the lowest level, all
computations are done with natural numbers in HOL Light. We improved
basic HOL Light procedures for natural numbers by representing natural
numerals over an arbitrary base (the base 2 is the standard base for
HOL Light natural numerals) and by providing arithmetic procedures for
computing with such numerals. Note that all computations required in
the nonlinear inequality verification procedure are done entirely
inside HOL Light, and the results of all arithmetic operations are HOL
Light theorems.  The native floating point operations of the computer
are not used in any formal proof.  As a consequence, the formal
verification of nonlinear inequalities in HOL Light is much slower
than the original informal C++ code.

The term \verb!the_nonlinear_inequalities! is defined as a conjunction
of several hundred nonlinear inequalities. The domains of these
inequalities have been partitioned to create more than 23,000
inequalities. The verification of all nonlinear inequalities in HOL
Light on the Microsoft Azure cloud took approximately 5000
processor-hours. Almost all verifications were made in parallel with
32 cores using GNU parallel \cite{Tange2011a}.  Hence the real time
was less than a week ($5000 < 32\times 168$ hours per week). These
verifications were made in July and August, 2014.  Nonlinear
inequalities were verified with compiled versions of HOL Light and the
verification tool developed in Solovyev's 2012 thesis.

The verifications were rechecked at Radboud University on 
60 hyperthreading Xeon 2.3GHz CPUs, in October 2014.  This
second verification required about 9370 processor-hours over a period of six
days.  Identical results were obtained in these repeated calculations.

\section{Combining HOL Light sessions}

The nonlinear inequalities were obtained in a number of separate sessions
of HOL Light that were run in parallel.  By the design of HOL Light, it
is not possible to pass a theorem from one session to another without
fully reconstructing the proof in each session.  
To combine the results into a single session of HOL Light, we used a
specially modified version of HOL Light that accepts a theorem from
another session without proof.  We briefly describe this modified
version of HOL Light.

Each theorem is expressed by means of a collection of constants, and
those constants are defined by other constants, recursively extending
back to the primitive constants in the HOL Light kernel.  Similarly,
the theorem and constants have types, and those types also extend
recursively back through other constants and types to the primitive
types and constants of the kernel.  A theorem relies on a list of
axioms, which also have histories of constants and types.  The
semantics of a theorem is determined by this entire history of
constants, types, and axioms, reaching back to the kernel.

The modified version of HOL Light is designed in such a way that a
theorem can be imported from another session, provided the theorem is
proved in another session, and the entire histories of constants,
types, and axioms for that theorem are exactly the same in the two
sessions.  To implement this in code, each theorem is transformed into
canonical form.  To export a theorem, the canonical form of the
theorem and its entire history are converted faithfully to a large
string, and the MD5 hash of the string is saved to disk.  The modified
version of HOL Light then allows the import of a theorem if the
appropriate MD5 is found.  The import mechanism is prevented from
being used in ways other than the intended use, through the scoping
rules of the OCaml language.

To check that no pieces were overlooked in the distribution of
inequalities to various cores, the pieces have been reassembled in the
specially modified version of HOL Light.  In that version, we obtain
a formal proof of the theorem

\begin{verbatim}
|- the_nonlinear_inequalities
\end{verbatim}

This theorem is exactly the assumption made in the formal proof of the
Kepler conjecture, as stated in Section~\ref{sec:statement}.  We
remark that the modified version of HOL Light is not used during the
proof of any other result of the Kepler conjecture.  It is only used
to assemble the small theorems from parallel sessions to produce this
one master theorem.

\section{Tame Classification}\label{sec:tc}


The first major success of the Flyspeck project was the formalization
of the classification of tame plane graphs. In the original proof of
the Kepler conjecture, this classification was done by computer, using
custom software to generate plane graphs satisfying given
properties. This formalization project thus involved the verification
of computer code. The formal verification of the code became the
subject of Gertrud Bauer's PhD thesis under the direction of
Nipkow. The work was completed by Nipkow \cite{NipkowBS-IJCAR06}.

As explained in Section~\ref{sec:tf}, the tame plane graphs encode
the possible counterexamples to the Kepler conjecture as plane
graphs. The {\it archive} is a computer-generated list of all tame
graphs. It is a text file that can be imported by the different parts
of the proof. In this section we explain how the following
completeness theorem is formalized in Isabelle/HOL:

\begin{lstlisting}[keepspaces=true,stringstyle=\tt,basicstyle=\small,%
frame=none,framesep=8pt,mathescape,morekeywords={and,shows},columns=flexible]
|- "g $\in$ PlaneGraphs" and "tame g" shows "fgraph g $\in_{\simeq}$ Archive"
\end{lstlisting}

The meaning of the
terms \verb!PlaneGraphs!, \verb!tame!, \verb!fgraph!, and
\verb!Archive! is explaned in the following paragraphs.
In informal terms, the completeness theorem 
asserts that every tame plane graph is
isomorphic to a graph appearing in the archive.  

Plane graphs are represented as an
$n$-tuple of data including a list of faces. Faces are represented
as lists of vertices, and each vertex is represented by an integer index.
A function \verb!fgraph! strips the $n$-tuple down to the list of faces.

To prove the completeness of the archive, we need to enumerate all
tame plane graphs. For this purpose we rely on the fact that HOL
contains a functional programming language. In essence, programs in
HOL are simply sets of recursion equations, i.e., pure logic.  The
original proof classifies tame plane graphs by a computer program
written in Java.  Therefore, as a first step in the formalization, we
recast the original Java program for the enumeration of tame plane graphs in
Isabelle/HOL. The result is a set of graphs called \verb!TameEnum!.  In
principle we could generate \verb!TameEnum! by formal proof but this
would be extremely time and space consuming because of the huge number
of graphs involved (see below). Therefore we rely on the ability of
Isabelle/HOL to execute closed HOL formulas by translating them
automatically into a functional programming language (in this case
ML), running the program, and accepting the original formula as a
theorem if the execution succeeds~\cite{DBLP:conf/flops/HaftmannN10}. The
programming language is merely used as a fast term rewriting engine.

We prove the completeness of the archive in two steps. First
we prove that every tame plane graph is in \verb!TameEnum!.
This is a long interactive proof in Isabelle/HOL.
Then we prove that every graph in \verb!TameEnum! is isomorphic to some
graph in the archive. Formally this can be expressed as follows:

\begin{lstlisting}[keepspaces=true,stringstyle=\tt,basicstyle=\small,%
frame=none,framesep=8pt,mathescape,morekeywords={and,shows},columns=flexible]
|- fgraph ` TameEnum $\subseteq_{\simeq}$ Archive
\end{lstlisting}
This is a closed formula that Isabelle proves automatically by evaluating it
(in ML) because all functions involved are executable.

In the following two subsections we give a high-level overview of the
formalization and proof. For more details see
\cite{NipkowBS-IJCAR06,Nipkow-ITP11}. The complete machine-checked proof,
including the archive is available online in the Archive of Formal Proofs
\url{afp.sf.net}~\cite{BauerN-AFP06}. Section~\ref{Archive} discusses
the size of the archive and some performance issues.

\subsection{plane graphs and their enumeration}\label{sec:tg}

Plane graphs are not defined by the conventional mathematical
definition. They are defined by an algorithm that starts with a
polygon and inductively adds loops to it in a way that intuitively
preserves planarity.  (The algorithm is an implementation in code of a
process of drawing a sequence of loops on a sheet of paper, each
connected to the previous without crossings.)

Expressed as a computer algorithm, the enumeration of plane graphs
proceeds inductively. It starts with a seed graph (the initial
polygon) with two faces (intuitively corresponding to the two
components in the plane of the complement of a Jordan curve), a final
outer one and a non-final inner one, where a {\it final} face means
that the algorithm is not permitted to make further modifications of
the face.  If a graph contains a non-final face, it can be subdivided
into a final face and any number of non-final ones.  Because a face
can be subdivided in many ways, this process defines a forest of
graphs. The leaves are final graphs.  The formalization defines an
executable function \verb!next_plane!  that maps a graph to the list
of successor graphs reachable by subdividing one non-final face.  The
set of plane graphs, denoted \verb!PlaneGraphs! in Isabelle/HOL, is
defined to be the set of final graphs reachable from some seed graph in
finitely many steps.

\subsection{tame graphs and their enumeration}

The definition of tameness is already relatively close to an
executable formulation. The two crucial constraints are that the faces
of a tame graph may only be triangles up to hexagons, and that the
``admissible weight'' of a tame graph is bounded from above. The
tameness constraints imply that there are only finitely many tame
plane graphs up to isomorphism (although we never need to prove this
directly).  The Isabelle/HOL predicate is called \verb!tame!.

The enumeration of tame plane graphs is a modified enumeration of plane
graphs where we remove final graphs that are definitely not tame, and prune
non-final graphs that cannot lead to any tame graphs.  The description
in~\cite{Hales-Annals} is deliberately sketchy, but the Java programs provide
precise pruning criteria. In the formalization we need to balance
effectiveness of pruning with simplicity of the completeness proof: weak
pruning criteria are easy to justify but lead to unacceptable run times of
the enumeration, sophisticated pruning criteria are difficult to justify
formally. In the end, for the formalization, we simplified the pruning criteria
found in the Java code.

The formalization defines a function \verb!next_tame! from a graph to
a list of graphs. It is a restricted version of \verb!next_plane!
where certain kinds of graphs are removed and pruned, as described
above.  For computational reasons, the tameness check here is
approximate: no tame graphs are removed, but non-tame graphs may be
produced. This is unproblematic: in the worst case a fake
counterexample to the Kepler conjecture is produced, which is
eliminated at a later stage of the proof, but we do not miss any real
ones.

Each step \verb!next_tame! is executable. The exhaustive enumeration
of all final graphs reachable from any seed graph via \verb!next_tame!
yields a set, denoted \verb!TameEnum!.

\subsection{the archive}
\label{Archive}

The archive that came with the original proof contained (for historical
reasons) over $5000$ graphs. The first
formalization~\cite{NipkowBS-IJCAR06} resulted in a reduced archive of
$2771$ graphs. During the completeness proof, the verified enumeration
of tame graphs has to go through $2\times 10^7$ intermediate and final
graphs, which takes a few hours. With the advent of the blueprint
proof, the definition of tameness changed. This change lead to $2
\times 10^9$ intermediate and final graphs and an archive with $18762$
graphs. The enumeration process had to be optimized to prevent it
running out space and time~\cite{Nipkow-ITP11}. In the end, the
enumeration again runs in a few hours.  The formalization uncovered
and fixed a bug in the original Java program, related to symmetry
optimizations.  Fortunately, this bug was not executed in the original
proof, so that the output was correct.  However, the bug caused two
graphs to be missed in an early draft of the blueprint proof.

\section{Importing results from Isabelle}

The tame graph classification was done in the Isabelle/HOL proof
assistant, while all the rest of the project has been carried out in
HOL Light.  It seems that it would be feasible to translate the
Isabelle code to HOL Light to have the entire project under the same
roof, but this lies beyond the scope of the Flyspeck project.

Current tools do not readily allow the automatic import of this result
from Isabelle to HOL Light.  A tool that automates the import from
Isabelle to HOL Light was written by McLaughlin with precisely this
application in mind~\cite{McLaughlin:2006:IJCAR}, but this tool has
not been maintained.  A more serious issue is that the proof in
Isabelle uses computational reflection as described at the end of
Section~\ref{sec:hl}, but the HOL Light kernel does not permit
reflection.  Thus, the reflected portions of the formal proof would
have to be modified as part of the import.

Instead, we leave the formalization of the Kepler conjecture
distributed between two different proof assistants.  In HOL Light, the
Isabelle work appears as an assumption, expressed through the
following definition.
\begin{obeylines}

\begin{verbatim}
|- import_tame_classification <=>
     (!g. g IN PlaneGraphs /\ tame g ==> fgraph g IN_simeq archive)
\end{verbatim}

\end{obeylines}
The left-hand side is exactly the assumption made in the formal proof
of the Kepler conjecture, as stated in Section~\ref{sec:statement}.
The right-hand side is the verbatim translation into HOL Light of the following
completeness theorem in Isabelle (repeated from above):




\begin{lstlisting}[keepspaces=true,stringstyle=\tt,basicstyle=\small,%
frame=none,framesep=8pt,mathescape,morekeywords={and,shows},columns=flexible]
|- "g $\in$ PlaneGraphs" and "tame g" shows "fgraph g $\in_{\simeq}$ Archive"
\end{lstlisting}


All of the HOL Light terms
\verb!PlaneGraphs!, \verb!tame!, \verb!archive!, \verb!IN_simeq!,
\verb!fgraph! are verbatim translations of the corresponding
definitions in Isabelle (extended recursively to the constants
appearing in the definitions).  The types are similarly translated
between proof assistants (lists to lists, natural numbers to natural
numbers, and so forth).  These definitions and types were translated
by hand.  The archive of graphs is generated from the same ML file for
both the HOL Light and the Isabelle statements.

Since the formal proof is distributed between two different systems
with two different logics, we briefly indicate why this theorem in
Isabelle must also be a theorem in HOL Light.  Briefly, this
particular statement could be expressed as a SAT problem in
first-order propositional logic.  SAT problems pass directly between
systems and are satisfiable in one system if and only if they are
satisfiable in the other (assuming the consistency of both systems).
In expressing the classification theorem as a SAT problem, the point
is that all quantifiers in the theorem run over bounded discrete sets,
allowing them to be be expanded as finitely many cases in
propositional logic.  We also note that the logics of Isabelle/HOL and
HOL Light are very closely related to each another.

\section{Linear programs}\label{sec:lp}

We return to the sketch of the proof from Section~\ref{sec:tf} and add
some details about the role of linear programming.  The blueprint
proof reduces infinite packings that are potential counterexamples to
the Kepler conjecture to finite packings $V$, called contravening
packings. The combinatorial structure of a contravening packing can be
encoded as a tame planar hypermap. The tame graph classification
theorem implies that there are finitely many tame planar hypermaps up
to isomorphism. For each such hypermap it is possible to generate a
list of inequalities which must be satisfied by a potential
counterexample associated with the given tame planar hypermap. Most of
these inequalities are nonlinear.

To rule out a potential counterexample, it is enough to show that the
system of nonlinear inequalities has no solution.  A linear relaxation
of these inequalities is obtained by replacing nonlinear quantities
with new variables. For instance, the dihedral angles of a simplex are
nonlinear functions of the edge lengths of the simplex; new variables
are introduced for each dihedral angle to linearize any expression that
is linear in the angles.  The next step is to show that
the linear program is not feasible. Infeasibility implies that the original
system of nonlinear inequalities is inconsistent, and hence there is
no contravening packing associated with the given tame planar
hypermap. The process of construction and solution of linear programs
is repeated for all tame planar hypermaps and it is shown that no
contravening packings (and hence no counterexamples to the Kepler
conjecture) exist.

There are two parts in the formal verification of linear
programs. First, linear programs are generated from formal theorems,
nonlinear inequalities, and tame planar hypermaps. Second, a general linear
program verification procedure verifies all generated linear programs.

The formal generation of linear programs follows the procedure
outlined above. The first step is generation of linear inequalities
from properties of contravening packings. Many such properties
directly follow from nonlinear inequalities. As described above,
nonlinear inequalities are transformed into linear inequalities by
introducing new variables for nonlinear expressions. In fact, we do
not change the original nonlinear inequalities but simply introduce
new HOL Light constants for nonlinear functions and reformulate
nonlinear inequalities in a linear form.

For example, suppose we have the following inequalities: $x +
x^2 \le 3$ and $x \ge 2$. Define $y = x^2$ to obtain the following
linear system of inequalities: $x + y \le 3$, $x \ge 2$, and $y \ge 4$.
(The last inequality follows from $x \ge 2$.)  We
ignore the nonlinear dependency of $y$ on $x$, and obtain a system of linear
inequalities which can be easily shown to be inconsistent.

The generation of linear inequalities from properties of contravening packings is
a semi-automatic procedure: many such inequalities are derived with a special
automatic procedure but some of them need manual formal proofs.

The next step is the relaxation of linear inequalities with irrational
coefficients. We do not solve linear programs with
irrational numbers directly.  Consider the example, $x - \sqrt{2} y \le
\pi$, $x + y \le \sqrt{35}$, $x \ge 5$, $y \ge 0$. These inequalities imply
\begin{equation*}\label{eqn:35}
x\ge 5, \quad
y \ge 0, \quad
x - 1.42 y \le 3.15,\quad
x + y \le 6.
\end{equation*}
This system with rational coefficients is inconsistent and thus
the original system is inconsistent.

The relaxation of irrational coefficients is done completely
automatically. In fact, inequalities with integer coefficients are
produced by multiplying each inequality by a sufficiently large power
of $10$ (decimal numerals are used in the verification of linear
programs).

The last step of the formal generation of linear programs is the
instantiation of free variables of linear inequalities with special
values computed from an associated tame planar hypermap. In this way, each
tame planar hypermap produces a linear program which is formally checked for
feasibility. Not all linear programs obtained in this way are
infeasible.  For about half of the tame planar hypermaps, linear relaxations
are insufficiently precise and do not yield infeasible linear
programs.  In that situation, a feasible linear program is split into
several cases where new linear inequalities are introduced. New cases
are produced by considering alternatives of the form $x \le a \vee a
\le x$ where $x$ is some variable and $a$ is a constant.
Case splitting leads to more precise linear relaxations.

Case splitting for verification of linear programs in the project is
automatic. In fact, a special informal procedure generates all required
cases first, and the formal verification procedure for linear programs is
applied to each case.

Formal verification of general linear programs is relatively easy. We
demonstrate our verification procedure with an example. A detailed
description of this procedure can be found in~\cite{Solovyev:LP}. All
variables in each verified linear program must be nonnegative and
bounded. In the following example, our goal is to verify that the following system is
infeasible:
\begin{equation*}
\begin{split}
x - 1.42 y \le 3.15,\ x + y \le 6,\ 
5 \le x \le 10,\ 0 \le y \le 10.
\end{split}
\end{equation*}
Introduce slack variables $s_1, s_2$ for inequalities which do not
define bounds of variables, and construct the following linear program:
\begin{equation}\label{eqn:lp}
\begin{split}
&\text{minimize $s_1 + s_2$ subject to}\\
&x - 1.42 y \le 3.15 + s_1,\ x + y \le 6 + s_2,\ 
5 \le x \le 10, 0 \le y \le 10,\ 
0 \le s_1,\ 0 \le  s_2.
\end{split}
\end{equation}
If the objective value of this linear program is positive then the
original system is infeasible. We are interested in a dual solution of
this linear program. Dual variables correspond to the constraints of
the primal linear program. We use an external linear programming tool
for finding dual solutions, which may be imprecise.  A procedure,
which is described in \cite{Solovyev:LP}, starts with an initial
imprecise dual solution and modifies it to obtain a carefully
constructed dual solution that has sufficient numerical precision to
prove the infeasibility of the original (primal) system.  In the
example at hand, we get the following modified dual solution:
\[
 (0.704, 1, -1.705, 0.001, -0.00032, 0, 0.296, 0),
\]
whose entries correspond with the ordered list of eight inequalities in the system
(\ref{eqn:lp}).  With this modified dual solution, the sum of the
constraints of the system (\ref{eqn:lp}) with the slack variables set
to zero, weighted by the coefficients of the dual vector, yields $0 x
+ 0 y \le -0.2974$. This contradiction shows that our original system
of inequalities is infeasible.  We stress that the coefficients of $x$
and $y$ in this sum are precisely zero, and this is a key feature of
the modified dual solutions.

If we know a modified dual solution, then the formal verification
reduces to the summation of inequalities with coefficients from the
modified dual solution and checking that the final result is
inconsistent. A modified dual solution can be found with informal
methods. We use GLPK for solving linear programs~\cite{website:GLPK}
and a special C\# program for finding the required modified dual
solutions for linear programs. All dual solutions are also converted
to integer solutions by multiplying all coefficients by a sufficiently
large power of $10$. The formal verification procedure uses formal
integer arithmetic. There are 43,078 linear programs (after
considering all possible cases). All these linear programs can be
verified in about $15$ hours on a 2.4GHz computer. The verification
time does not include time for generating modified dual
solutions. These dual solutions need only be computed once and saved
to files that are later loaded by the formal verification procedure.

\section{Auditing a distributed formal proof}

A proof assistant largely cuts the mathematical referees out of the
verification process.  This is not to say that human oversight is no
longer needed.  Rather, the nature of the oversight is such that
specialized mathematical expertise is only needed for a small
part of the process.  The rest of the audit of a formal proof can be
performed by any trained user of the HOL Light and Isabelle proof
assistants.

Adams \cite{adams2014flyspecking} describes the steps involved in the
auditing of this formal proof.  The proofs scripts must be executed to
see that they produce the claimed theorem as output.  The definitions
must be examined to see that the meaning of the final theorem (the
Kepler conjecture) agrees with the common understanding of the
theorem.  In other words, did the right theorem get formalized?  Were
any unapproved axioms added to the system?  The formal proof system
itself should be audited to make sure there is no foul play in the
syntax, visual display, and underlying internals.  A fradulent user of
a proof assistant might ``exploit a flaw to get the project completed
on time or on budget.  In their review, the auditor must assume
malicious intent, rather than use arguments about the improbability of
innocent error'' \cite{adams2014flyspecking}.

This particular formal proof has several special features that call
for careful auditing.  The most serious issue is that the full formal
proof was not obtained in a single session of HOL Light.  An audit
should check that the statement of the tame classification theorem in
Isabelle has been faithfully translated into HOL Light.  (It seems to
us that our single greatest vulnerability to error lies in the hand
translation of this one statement from Isabelle to HOL Light, but even
here there is no mathematical reasoning involved beyond a rote
translation.)  In particular, an audit should check that the long list
of tame graphs that is used in Isabelle is identical to the list that
is used in HOL Light.  (As mentioned above, both systems generated
their list from the same master file.)

The auditor should also check the
design of the special modification of HOL Light that was used to combine
nonlinear inequalities into a single session.

\section{Related work and acknowledgements}

Numerous other research projects in formal proofs have made use of
the Flyspeck project in some way or have been inspired by the needs of
Flyspeck.  These projects include the work cited above on linear programming
and formal proofs of nonlinear inequalities, automated translation of proofs
between formal proof systems~\cite{obua:import}, \cite{KaliszykK13}, 
\cite{McLaughlin:2006:IJCAR},
the refactoring of formal proofs \cite{adams2012recording}, machine learning
applied to proofs and proof automation ~\cite{KU14} \cite{KaliszykU2014}, an SSReflect mode for HOL Light~\cite{Solovyev-thesis}, and a mechanism to
execute trustworthy external arithmetic from
Isabelle/HOL~\cite{obua:phd},~\cite{ObuaN}.

The roles of the
various members of the Flyspeck project have been spelled out in the email
announcement of the project completion, posted
at~\cite{website:FlyspeckProject}. 

We wish to acknowledge the help, support, influence, and various
contributions of the following individuals:
Nguyen Duc Thinh,  
Nguyen Duc Tam, 
Vu Quang Thanh,
Vuong Anh Quyen,
%
Catalin Anghel, 
Jeremy Avigad, 
Henk Barendregt,
%
Herman Geuvers,
Georges Gonthier,
Daron Green,
Mary Johnston,
Christian Marchal,
Laurel 
Martin, 
%
Robert Solovay,
Erin Susick,
Dan Synek,
Nicholas Volker, 
Matthew Wampler-Doty, 
Benjamin Werner,
Freek Wiedijk, 
Carl Witty, and
Wenming Ye.

We wish to thank the following sources of institutional support: NSF
grant 0503447 on the ``Formal Foundations of Discrete Geometry" and NSF
grant 0804189 on the ``Formal Proof of the Kepler Conjecture,"
Microsoft Azure Research, William Benter Foundation, University of
Pittsburgh, Radboud Research Facilities, Institute of Math (VAST), and
VIASM.

\section{Appendix on definitions}\label{sec:ap}

The following theorem provides evidence that key definitions in the
statement of the Kepler conjecture are the expected ones.

\begin{obeylines}

\begin{verbatim}
|-  
// real absolute value:
   (&0 <= x ==> abs x = x) /\ (x < &0 ==> abs x = -- x) /\   

// powers:
    x pow 0 = &1 /\ x pow SUC n = x * x pow n /\

// square root:
   (&0 <= x ==> &0 <= sqrt x /\ sqrt x pow 2 = x) /\ 

// finite sums:
   sum (0..0) f = f 0 /\ sum (0..SUC n) f =  
     sum (0..n) f + f (SUC n) /\ 

// pi:
   abs (pi / &4 - sum (0..n) (\i. (-- &1) pow i / &(2 * i + 1))) 
     <= &1 / &(2 * n + 3) /\

// finite sets and their cardinalities:
   (A HAS_SIZE n <=> FINITE A /\ CARD A = n) /\
   {} HAS_SIZE 0 /\ {a} HAS_SIZE 1 /\ 
   (A HAS_SIZE m /\ B HAS_SIZE n /\ (A INTER B) HAS_SIZE p 
     ==> (A UNION B) HAS_SIZE (m+n - p)) /\

// bijection between R^3 and ordered triples of reals:
   triple_of_real3 o r3 = (\w:real#real#real. w) /\
   r3 o triple_of_real3 = (\v:real^3. v) /\ 

// the origin:
   vec 0 = r3(&0,&0,&0) /\

// the metric on R^3:
   dist(r3(x,y,z),r3(x',y',z')) = 
     sqrt((x - x') pow 2 + (y - y') pow 2 + (z - z') pow 2) /\

// a packing:
   (packing V <=> 
     (!u v. u IN V /\ v IN V /\ dist(u,v) < &2 ==> u = v))
\end{verbatim}

\end{obeylines}

      \bibliographystyle{plainnat}
\bibliography{announce_refs}




\end{document}